\documentclass[12pt]{article}
\usepackage{amssymb,amsfonts}
\usepackage{graphicx}
\usepackage{amscd}

\newtheorem{lemma}{Lemma}
\newtheorem{prop}{Proposition}

\title{Rolling Balls and Octonions}
\author{A.~A.~Agrachev\thanks{SISSA--ISAS, Trieste \& Steklov Math. Inst.,
Moscow}}
\date{}
\begin{document}
\maketitle
\begin{abstract}
In this semi-expository paper we disclose hidden symmetries of a
classical nonholonomic kinematic model and try to explain
geometric meaning of basic invariants of vector distributions.
\end{abstract}

\section{Introduction}
The paper is written to the 70th anniversary of Vladimir Igorevich
Arnold. This is just a small mathematical souvenir but I hope that
Vladimir Igorevich will get some pleasure looking over it. The
content of the paper is well described by the cryptogram below.
The figure represents the roots system of the exceptional group
$G_2$ (the automorphisms group of octonions) and two touching each
other circles whose diameters are in ratio 3:1.
\begin{center}
\includegraphics[scale=0.25]{logo.1}
\end{center}

Our starting point is a classical nonholonomic kinematic system
that is rather important in robotics: a rigid body rolling over a
surface without slipping or twisting. The surface is supposed to
be the surface of another rigid body so that the situation is, in
fact, symmetric: one body is rolling over another one. We also
assume that the surfaces of the bodies are smooth and cannot touch
each other in more than one point.

This system has a 5-dimensional configuration space: coordinates
of the points on each surface where the surfaces touch each other
give 4 parameters; fifth parameter measures mutual orientation of
the bodies at the touching point.

Now let us imagine one of the bodies immovable and another
rolling. Given an initial configuration, it is rather clear that
one can roll the body in a unique way along any curve on the
immovable surface starting from the initial touching point. In
other words, given an initial configuration, admissible motions
are parameterized by the curves on the two-dimensional surface
with a fixed initial point. On the other hand, it is not hard to
prove (see \cite[Ch.\,24]{as}) that admissible motions allow to
reach any configuration from any other as soon as two bodies are
not congruent. We are thus in a typical nonholonomic situation
since the configuration space is 5-dimensional.

Now turn to Mathematics. Admissible velocities (i.e. velocities of
admissible motions) form a rank 2 vector distribution on the
configuration space. This distribution is involutive if and only
if our bodies are the balls of equal radii. More precisely, let
vector fields $f,g$ form a local basis of our distribution,
$[f,g]$ be the commutator (Lie bracket) of the fields, and $q$ a
point of the configuration space. It happens that vectors $f(q),
g(q), [f,g](q)$ are linearly independent if and only if the
curvatures of two surfaces are not equal at their touching points
corresponding to the configuration $q$. Moreover, if these
curvatures are not equal, then:
$$
f(q)\wedge g(q)\wedge [f,g](q)\wedge [f,[f,g]](q)\wedge
[g,[g,f]](q) \ne 0. \eqno (1)
$$
In other words, the basic fields, their first and second order Lie
brackets are all linearly independent.

The germs of rank 2 distributions in $\mathbb R^5$ with property
(1) were first studied by Elie Cartan in his famous paper
\cite{ca}. A brief count of parameters demonstrates that 5 is the
lowest dimension where the classification of generic germs of
distributions must have functional invariants. This happens for
rank 2 and rank 3 distributions in $\mathbb R^5$; moreover, the
classifications for both values of the rank are essentially
equivalent and it is sufficient to study the case of rank 2.

Cartan found a fundamental tensor invariant that is a degree 4
symmetric form on the distribution: the desired functional
invariant is just the cross-ratio of the roots of this form. There
is exactly one equivalence class for which the Cartan's form is
identical zero. We call the germs from this class {\it flat}; a
germ of the distribution is flat if and only if it admits a basis
generating the 5-dimensional nilpotent Lie algebra. Cartan showed
that symmetry group of a flat distribution is the split form of
the 14-dimensional exceptional simple Lie group $G_2$. In all
other cases, dimension of the symmetry group does not exceed 7.

How all that is related to the rolling bodies? It is rather
obvious that  the group of symmetries of the ``no slipping no
twisting distribution" is trivial for generic bodies. This group
acts transitively on the configuration space if and only if the
surfaces of both bodies have constant curvatures. Let us stay with
the case of constant nonnegative curvatures so that the surfaces
are spheres (one of them may be a plane, i.e. the sphere of the
infinite radius).

Natural symmetries are isometries of the spheres; they form a
6-dimen\-si\-o\-nal Lie group. It was Robert Bryant who claimed
first that no slipping no twisting distribution is flat in the
case of spheres whose radii are in ratio 3:1. He insisted that
just followed the Cartan's method and never published this fact in
a paper or a preprint. It also remained absolutely unclear what
are hidden additional symmetries (dimension of the symmetry group
jumps from 6 to 14 when ratio of the radii gets 3:1!).
Unfortunately, Cartan's method does not give much in this regard:
the construction of the fundamental invariant is based on the
involved reduction--prolongation procedures in the jet spaces, and
any connection with the original geometric problem is lost quite
far from the end of the way.

The goal of this note is to finally untwine this puzzle. In
Section~2 we give a simple twistor model for the configuration
space of the rolling balls problem and for the no slipping no
twisting distribution; the role of the group $G_2$ is not yet
clear from this model.

Split form of $G_2$ is the group of automorphisms of {\it
split-octonions}, a hyperbolic version of usual octonions, where
positive definite quadratic form (square of the norm of octonions)
is substituted by a non-degenerate sign-indefinite quadratic form.
Split-octonions have nontrivial divisors of zero (zero locus of
the sign-indefinite form). A simple quadratic transformation
demonstrates that no slipping no twisting distribution is
equivalent to the ``divisors of zero distribution" of
split-octonions if the radii of the balls are in ratio 3:1. This
is the subject of Section~3.

In Section~4 we outline a recently developed variational way to
construct differential invariants of vector distributions in order
to put the rolling bodies model in a broader framework and to
explain the meaning of basic invariants.

\medskip\noindent \textsl{Acknowledgements.} I am grateful to Robert Bryant,
Aroldo Kaplan and Igor Zelenko for very interesting stimulating
discussions.

\section{Twistor model}

We study admissible motions of two balls $B_r$ and $B_R$ rolling
one over another without slipping or twisting. Here $r$ and $R$
are the radii of the balls. Instant configuration of the system of
two balls is determined by an orientation preserving isometry of
the tangent planes to the spheres $S_r=\partial B_r$ and
$S_R=\partial B_R$ at the points where the balls touch each other.
In other words, the state space of our kinematic system is:
$$
M_{R,r}=\left\{\mu:T_{q_1}S_r\to T_{q_2}S_R \mid {q_1\in
S_r,\,q_2\in S_R,\ \mu\ \mathrm{is\ an}\atop \mathrm{isometry\ of\
oriented\ planes}}\right\}.
$$

It is easy to see that $M_{R,r}$ is a smooth 5-dimensional
manifold. Motions of the system are families of isometries
$\mu(t):T_{q_1(t)}S_r\to T_{q_2(t)}S_R,\ t\in\mathbb R$. The no
slipping condition reads:
$$
\mu(t)(\dot q_1(t))=\dot q_2(t).
$$
The no twisting condition requires that $\mu(t)$ transforms
parallel along $q_1(t)$ vector fields in the parallel along
$q_2(t)$ vector fields.

These two conditions define a rank~2 vector distribution $D^{R,r}$
on $M_{R,r}$. We have: $D^{R,r}=\bigcup\limits_{\mu\in
M_{R,r}}D^{R,r}_\mu$, where $D^{R,r}_\mu$ is a two-dimensional
subspace of $T_\mu M_{R,r}$; admissible motions of the two balls
system are exactly integral curves of the distribution $D^{R,r}$.
Given an initial configuration, the ball $B_r$ can be rolled in a
unique way along any smooth curve on $S_R$ and the same is true if
we transpose $r$ and $R$. In the formal geometric language this
observation just means that the subspace $D^{R,r}_\mu\subset T_\mu
M_{R,r},$ where $\mu:T_{q_1}\to T_{q_2}$, is projected one-to-one
onto $T_{q_1}S_R$ and $T_{q_2}S_r$, $\forall \mu\in M_{R,r}$.

In what follows, we treat $S_r$ and $S_R$ as unit spheres in
$\mathbb R^3$ with re-scaled metrics:
$$
M_{R,r}=\{\mu:q_1^\bot\to q_2^\bot \mid q_i\in\mathbb R^3,\
|q_i|=1,\ i=1,2,\ R|\mu(v)|=r|v|,\ \forall v\in q_1^\bot\}.
$$
Let $\rho=\frac Rr$, the homothety $\iota_\rho:\mu\mapsto\rho\mu,\
\mu\in M_{R,r}$, transforms $M_{R,r}$ in $M_{1,1}$. We set
$D^\rho=\iota_{\rho *}D^{R,r}$. The distribution $D^\rho$ on
$M_{1,1}$ is determined by the ``re-scaled no slipping condition"
$$
\mu(\dot q_1(t))=\rho\dot q_2(t)
$$
and the no twisting condition; the last one remains unchanged.

From now on we will deal with fixed space $M_{1,1}$ endowed with
the family of distributions $D^\rho$ instead of the family of
pairs $\left(M_{R,r},D^{R,r}\right)$. In order to explicitly
describe the distributions $D^\rho$ we use a classical
quaternion's parameterization of spherical bundle $\mathfrak
p:\mathcal S\to S^2$, to the unit sphere $S^2\subset\mathbb R^3$,
where
$$
\mathcal S=\{(q,v)\in\mathbb R^3\times\mathbb R^3 \mid |q|=|v|=1,\
\langle q,v\rangle=0\},\quad \mathfrak p(q,v)=q.
$$
Let us recall this parameterization. We identify $\mathbb R^3$
with the space of imaginary quaternions:
$$
\mathbb R^3=\{\alpha i+\beta j+\gamma k\mid
\alpha,\beta,\gamma\in\mathbb R\}\subset\mathbb H.
$$
Let $S^3=\{w\in\mathbb H\mid |w|=1\}$ be the group of unitary
quaternions; then $\mathfrak h:S^3\to S^2,\ \mathfrak h(w)=\bar
wiw$ is the classical Hopf bundle, while the mapping
$$
\Psi:S^3\to\mathcal S,\quad \Psi(w)=(\bar wiw,\bar wjw)
$$
is a double covering. Moreover, the diagram
$$
\begin{CD}  S^3 @>\Psi>>\mathcal S\\
@VV\mathfrak{h}V @VV\mathfrak{p}V\\
S^2 @= S^2
\end{CD}
$$
is commutative, hence $\Psi$ is a fiber-wise mapping of the bundle
$\mathfrak h:S^3\stackrel{S^1}{\longrightarrow}S^2$ onto the
bundle $\mathfrak p:\mathcal S\stackrel{S^1}{\longrightarrow}S^2$
which induces a double covering of the fibers. The fibers of the
bundle $\mathfrak h:S^3\to S^2$ are the residue classes
$\{e^{i\theta}w\mid \theta\in\mathbb R\,\mathrm{mod}\,2\pi\}$ of
the one-parametric subgroup generated by $i$.

The normal to the fibers distribution
$$
span\,\{jw,kw\}\subset T_wS^3,\quad w\in S^3
$$
is a connection on the principal bundle $\mathfrak h:S^3\to S^2$.
It is easy to see that $\Psi$ transforms this connection in the Levi
Civita connection on the bundle $\mathfrak p:\mathcal S \to S^2$,
which defines standard parallel transport on $S^2$.

We are now ready to give a quaternion's model of the rolling balls
configuration space $M_{1,1}$ (more precisely, of a double
covering to $M_{1,1}$) and of the no slipping no twisting
distributions $D^\rho$. For any $w_1,w_2\in S^3$ there exists a
unique orientation preserving isometry of the fiber
$$
\mathfrak h^{-1}(\mathfrak
h(w_1))=\{e^{i\theta}w_1\mid\theta\in\mathbb
R\,\mathrm{mod}\,2\pi\}
$$
on the fiber $\mathfrak h^{-1}(\mathfrak h(w_2))$ that sends $w_1$
to $w_2$. This isometry sends $e^{i\theta}w_1$ to
$e^{i\theta}w_2$. Moreover, pairs $(w_1,w_2)$ and $(w_1',w_2')$
define one and the same isometry if and only if
$w_1'=e^{i{\theta}'}w_1,\ w_2'=e^{i{\theta}'}w_1$ (with one and
the same ${\theta}'$). Hence the coset space of $S^3\times S^3$ by
the action $(w_1,w_2)\mapsto (e^{i\theta}w_1,e^{i\theta}w_2)$ of
the one-parametric group $\{e^{i\theta}\mid\theta\in\mathbb
R\,\mathrm{mod}\,2\pi\}$ is a double covering of $M_{1,1}$. We use
symbol $\mathbf M$ for this coset space and $\pi:S^3\times S^3\to
\mathbf M$ for the canonical projection. By $\mathbf
D^\rho=\bigcup\limits_{\mathbf x\in\mathbf D}\mathbf
D^\rho_{\mathbf x}$ we denote the rank 2 distribution on $\mathbf
M$ that is the pullback of $D^\rho$ by the double covering
$\mathbf M\to M_{1,1}$. Then
$$
\mathbf D^\rho_{\pi(w_1,w_2)}=\pi_*span\,\{(jw_1,\rho
jw_2),\,(kw_1,\rho kw_2)\}\quad \forall w_1,w_2\in S^3.
$$

\bigskip Let us now treat the quaternionic space $\mathbb
H^2=\{(w_1,w_2) \mid w_i\in\mathbb H\}$ as $\mathbb C^4$, where
$w_1=z_1+z_2j,\ w_2=z_3+z_4j,\ z_l\in\mathbb C,\ l=1,\ldots,4$. We
see that $\mathbf M$ is nothing else but a complex projective
conic,
$$
\mathbf M=\{z_1:z_2:z_3:z_4 \mid |z_1|^2+|z_2|^2=|z_3|^2+|z_4|^2\}
\subset\mathbb{CP}^3.
$$ This conic is often called ``the space of isotropic twistors".
Moreover,
$$
\mathbf D^\rho_{\pi(w_1,w_2)}=\pi_*\mathbb Cj(w_1,\rho
w_2).
$$

\section{Split-Octonions}

We also treat $\mathbb H^2$ as the algebra $\hat\mathbb
O=\{w_1+\ell w_2 \mid w_i\in\mathbb H\}$ of split-octonions, where
$$
(a+\ell b)(c+\ell d)=(ac+d\bar b)+\ell(\bar ad+cb). \eqno (2)
$$ Let
$x=w_1+\ell w_2$, $\bar x=\bar w_1-\ell\bar w_2$ and $Q(x)=\bar
xx=|w_1|^2-|w_2|^2$. Then $Q(xy)=Q(x)Q(y)$ and $x^{-1}=\frac{\bar
x}{Q(x)}$ as soon as $Q(x)\ne 0$; the cone $Q^{-1}(0)$ consists of
the divisors of zero.

The automorphisms group of the algebra $\hat\mathbb O$ is the split-form
of the exceptional Lie group $G_2$. These automorphisms preserve quadratic form
$Q$ and hence its polarization $\mathbf Q(x,y)=\frac 14(Q(x+y)-Q(x-y))$. In particular,
these automorphisms preserve the subspace
$\mathbb R^7=\{x\in\hat\mathbb O \mid \mathbf Q(1,x)=0\}$ and the conic
$$
K=\{x\in\mathbb R^7 \mid Q(x)=0\}=\{x\in\hat\mathbb O \mid xx=0\}.
$$
Moreover, the
automorphisms group of $\hat\mathbb O$ acts transitively on the
``spherization" $\mathbf K=\{\mathbb R_+x \mid x\in K\setminus 0
\}=S^2\times S^3$ of the cone $K$.

To any $x\in K\setminus 0$ there associated a 3-dimensional
subspace of divisors of 0 $$\Delta_x=\{y\in\mathbb R^7 \mid
xy=0\};$$ the ``spherization" turns $\Delta_x$ into a
2-dimensional subspace $\mathbf\Delta_{\mathbf x}\subset
T_{\mathbf x}\mathbf K$, where $\mathbf x=\mathbb R_+x$.
Obviously, the automorphisms group of $\hat\mathbb O$ preserves
vector distribution $\mathbf \Delta=\{\mathbf\Delta_{\mathbf
x}\}_{\mathbf x\in\mathbf K}$.

\medskip
\begin{prop} The mapping
$\Phi:(w_1+\ell w_2)\mapsto(w_1^{-1}iw_1+\ell(w_1^{-1}w_2))$
induces the diffeomorphism of $\mathbf M$ onto $\mathbf K$.
Moreover, the differential of this diffeomorphism transforms the
``no slipping, no twisting" distribution $\mathbf D^3$ in the
``divisors of zero" distribution $\mathbf\Delta$.
\end{prop}
{\bf Proof.} Let $\hat\Phi:\mathbf M\to\mathbf K$ be the mapping
induced by $\Phi$. We give an explicit formula for
$\hat\Phi^{-1}$: take $v_1\in S^2,\ v_1=\mathfrak h(w_1)$, then
$$
\hat\Phi^{-1}(v_1+\ell v_2)=\pi(w_1+\ell(w_1v_2)),\quad \forall
v_2\in S^3.
$$

Now we have to prove that $\Phi(x)\left(D_x\Phi y\right)=0$
for any $x=w_1+\ell w_2,\ y=zjw_1+3\ell(zjw_2)$
such that $|w_1|=|w_2|,\ z\in\mathbb
C$. It is sufficient to make calculation in the case
$|w_1|=|w_2|=1$. We have:
$$
D_{w_1+\ell w_2}\Phi\left(zjw_1+\ell(zjw_2)\right)=2\bar
w_1zkw_1+2\ell(\bar w_1zjw_2).
$$
The desired result now follows from the multiplication rule (2).

\bigskip Let us give an explicit parameterization of the distribution
$\mathbf \Delta$ on $\mathbf K$. First we parameterize $\mathbf K$ itself:
$$
\mathbf K=\{v_1+\ell v_2\mid v_1\in\mathbb R^3,\,v_2\in\mathbb H,\
|v_1|=|v_2|=1\}.
$$
Then $\mathbf\Delta_{v_1+\ell v_2}=\{v_1u+\ell(uv_2)\mid
u,(v_1u)\in\mathbb R^3\}$; this is a simple corollary of the
multiplication rule (2). Let $v\in\mathbb R^3,\ |v|=1$; the
mapping $w\mapsto w+vwv$ maps $\mathbb H$ onto the subspace
$\{u\mid u,(vu)\in\mathbb R^3\}$. Now we substitute $u$ by
$w+v_1wv_1$ in the above description of $\mathbf\Delta_{v_1+\ell
v_2}$ and obtain the final parameterization:
$$
\mathbf\Delta_{v_1+\ell v_2}=\{[v_1,w]+\ell((w+v_1wv_1)v_2)\mid w\in\mathbb H\}.
$$

\section{Jacobi curves}

In this section we briefly describe the variational approach to
differential invariants of vector distributions
(see \cite{aII,az,ze}) in order to put the rolling balls model in a wider
perspective. This approach is based on the contemporary optimal
control techniques and suggests an alternative to the classical
equivalence method (see \cite{ca,ga,mo}) of Elie Cartan.

Rank $k$ vector distribution $\Delta$ on the $n$-dimen\-si\-o\-nal
smooth manifold $M$ is just a smooth vector subbundle of the tangent
bundle $TM$:
$$
\Delta=\bigcup\limits_{q\in M}\Delta_q,\quad \Delta_q\subset T_qM,\quad
\dim\Delta_q=k.
$$
Distributions $\Delta$ and $\Delta'$ are called locally equivalent
at $q_0\in M$ if there exists a neighborhood $O_{q_0}\subset M$ of
$q_0$ and a diffeomorphism $\Phi:O_{q_0}\to O_{q_0}$ such that
$\Phi_*\Delta_q=\Delta'_{\Phi(q)},\ \forall q\in O_{q_0}$.

A local basis of $\Delta$ is a $k$-tuple of smooth vector fields
$f_1,\ldots,f_k\in\mathrm{Vec}\,M$ such that
$$
\Delta_q=span\{f_1(q),\ldots,f_k(q)\},\ q\in O_{q_0}.
$$
Given a local basis, one may compute the {\it flag of the
distribution}:
$$
\Delta^l_q=span\left\{(\mbox{ad}f_{i_j}\cdots\mbox{ad}f_{i_1}f_{i_0})(q) \mid
0\le j<l\right\},\ l=1,2,\ldots,
$$
where $\mbox{ad}f\,g\stackrel{def}{=}[f,g]$
is the Lie bracket.

It is easy to see that the subspaces
$\Delta_q^l$ do not depend on the local basis. We set
$\Delta^l=\bigcup\limits_{q\in M}\Delta^l_q$, a growing sequence of subsets in
$TM$. This sequence stabilizes as soon as $\Delta^{l+1}=\Delta^l$. The
distribution is involutive if and only if $\Delta^2=\Delta$, it is completely
nonholonomic (our subject) if $\Delta^l=TM$ for sufficiently big $l$. Generic
distributions are, of course, completely nonholonomic.

Integral curves of the distribution are often called {\it horizontal paths}.
It is convenient to consider all paths of class $H^1$,
not only smooth ones. We thus have a Hilbert manifold $\Omega_{\Delta}$
of horizontal paths:
$$
\Omega_\Delta=\{\gamma\in H^1([0,1];M)
\mid \dot\gamma(t)\in\Delta_{\gamma(t)},\mathrm{for\ a.e.}\ t\in[0,1]\}.
$$
Now consider {\it boundary mappings}
$$\partial_t:\Omega_\Delta\to M\times M$$ defined by the formula
$
\partial_t:(\gamma)=(\gamma(0),\gamma(t)).$
It is easy to show that $\partial_t$ are smooth mappings.

Critical points of the mapping $\partial_1$ are called {\it
singular curves} of $\Delta$. Any singular curve is automatically a critical
point of $\partial_t,\ \forall t\in[0,1].$
Moreover, any singular curve possesses a {\it singular extremal},
i.\,e. an $H^1$-curve $\lambda:[0,1]\to T^*M$ in the cotangent bundle to
$M$ such that
$$\lambda(t)\in T^*_{\gamma(t)}M\setminus\{0\},\quad
(\lambda(t),-\lambda(0))D_\gamma
\partial_t=0,\quad \forall t\in[0,1].
$$
We set:
$$\Delta_q^\bot=\{\nu\in T_q^*M \mid \langle\nu,\Delta_q\rangle=0,\nu\ne0\},\quad
\Delta^\bot=\bigcup\limits_{q\in M}\Delta_q^\bot.$$
Obviously, $\Delta^\bot$ is a smooth $(n+k)$-dimensional submanifold of $T^*M$
(the annihilator of $\Delta$).

Let $\sigma$ be the canonical symplectic structure on $T^*M$.
Pontryagin Maximum Principle implies that
a curve $\lambda$ in $T^*M$ is a singular extremal
if and only if it is a characteristics of the form
$\sigma\bigr|_{\Delta^\bot}$; in other words,
$$ \lambda(t)\in\Delta^\bot,\quad
\dot\lambda(t)\in\ker\left(\sigma\bigr|_{\Delta^\bot}\right),\quad
0\le t\le 1.
$$

All singular extremals are contained in the {\it characteristic variety}
$$
C_\Delta=\left\{z\in\Delta^\bot \mid
\ker\sigma_z\bigr|_{\Delta^\bot}\ne 0\right\}.
$$
We have: $C_\Delta=\Delta^{2\bot}$ if $k=2$;
$C_\Delta=\Delta^\bot$ if $k$ is odd; typically, $C_\Delta$ is a
codim 1 submanifold of $\Delta^\bot$ if $k$ is even.

A complete description of singular extremals is a hard task;
to simplify the job we focus only on the
regular part of the characteristic variety. We set:
$$
C^0_\Delta=\left\{z\in
C_\Delta \mid \dim\,\ker\sigma_z\bigr|_{\Delta^\bot}\le 2,\
\dim\,\ker\sigma_z\bigr|_{\Delta^\bot}\cap
T_zC_\Delta=1\right\}.
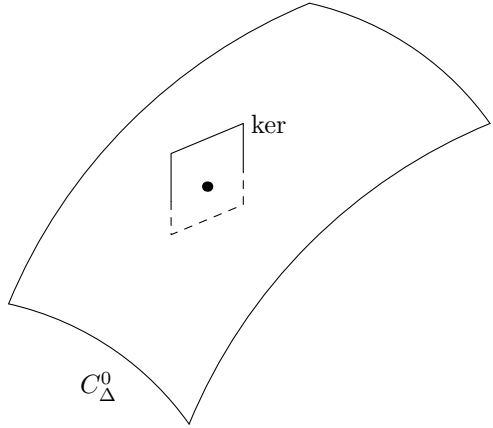
$$
\begin{center}
\includegraphics{charact.2}
\end{center}
If $k=2$, then $C^0_\Delta=\Delta^{2\bot}\setminus\Delta^{3\bot}$.

$C^0_\Delta$ is a smooth submanifold of $\Delta^\bot$;
it is foliated by singular extremals and by
the fibers $T^*_qM\cap C^0_\Delta$.

The motion along singular extremals defines a local flow on
$C^0_\Delta$; typically, this flow is not fiber-wise, i.\,e. it
does not transform fibers into fibers (see the figure above).

Given $z\in C^0_\Delta$ and an appropriate small neighborhood
$\mathcal C^0_z$  of $z$ in $C^0_\Delta$ we consider the
canonical projection
$$
F:\mathcal C_z^0\longrightarrow\mathcal C^0_z/\{\mathrm{singular\
extremals\ foliation}\}
$$
of $\mathcal C_z^0$ on the space of contained in $\mathcal C_z^0$
singular extremals.

Assume that $\lambda$ is a singular extremal through $z$ and it is
associated to a singular curve $\gamma$, i.\,e. $\lambda(0)=z,\
\lambda(t)\in T^*_{\gamma(t)}M$. Consider a family of subspaces
$$
J^0_\lambda(t)=T_\lambda F(T^*_{\gamma(t)}M\cap\mathcal C^0_z)
$$
of the space
$$
T_\lambda\mathcal C^0_z/\{\mathrm{singular\ extremals\ foliation}\}\cong
T_zC^0_\Delta/T_z\lambda.
$$
Then $t\mapsto J^0_\lambda(t)$ is a smooth curve in the
correspondent Grassmann manifold.
Geometry of the curves $t\mapsto J^0_\lambda(\cdot)$
reflects the dynamics of the fibers along
singular extremals and contains the fundamental information about
distribution $\Delta$.

\bigskip
In what follows we assume that $k=2,\ n\ge 5$ and $\Delta_q^2,\
\Delta_q^3$ have maximal possible dimensions, i.\,e.
$\dim\Delta_q^2=3,\ \dim\Delta_q^3=5$.

\medskip\noindent {\bf i.}
First we consider the case $n=5$ that is the case studied by
Cartan (see Introduction). Let $z\in T^*_qM$ and $\pi:T_z(T^*M)\to
T_qM$ be the differential at $z$ of the projection $T^*M\to M$;
then $\pi(J^0_\lambda(t))\subset z^\bot\subset T_qM$. Moreover,
$T_q\gamma\subset\pi(J^0_\lambda(t))$ and
$t\mapsto\pi\left(J^0_\lambda(t)\right)$ is a curve in the
projective plane $\mathbb P\left(z^\bot/T_q\gamma\right)$.

\begin{prop}[see \cite{aznu}]
A rank 2 distribution $\Delta$ on the 5-dimensional manifold is
flat if and only if the curve $\pi\left(J^0_\lambda(\cdot)\right)$
is a quadric, for any singular extremal $\lambda$. .
\end{prop}

In general, $\pi\left(J^0_\lambda(\cdot)\right)$ is not a quadric.
Let $K_z(q)\subset z^\bot$ be the best approximating quadric for
the curve $\pi\left(J^0_\lambda(\cdot)\right)$ near the point
corresponding to the zero value of parameter $t$ (the {\it
osculating quadric} of classical projective geometry); then
$K_z(q)$ is zero locus of a signature $(2,1)$ quadratic form on
$z^\bot/T_q\gamma$. We can, of course, treat $K_z(q)$ as zero
locus of a degenerate quadratic form on $z^\bot$. Finally,
$\mathcal K(q)=\bigcup\limits_{z\in\Delta_q^{2\bot}}K_z(q)$ is
zero locus of a $(3,2)$ quadratic form on $T_qM$ (see \cite{aznu}
for details).

The family of quadratic cones
$\mathcal K(q),\ q\in M$ is an intrinsically ``raised" from
$\Delta$ conformal structure on $M$ and $\Delta_q\subset K(q)$.
This conformal structure was first found by Nurowski \cite{nu} who used the
Cartan's equivalence method.

\medskip\noindent {\bf Remark.} Nurowski's conformal structure has
a particularly simple description for the divisors of zero distribution
$\mathbf\Delta$ of Sec.~3. Namely,
$\mathcal K(\mathbf x)=Q^{-1}(0)\cap T_{\mathbf x}\mathbf K,\
\mathbf x\in\mathbf K$, in this case.

\bigskip\noindent {\bf ii.}
From now on $n$ is any integer greater or equal to 5.
Let $z\in C^0_\Delta$, $\lambda$ be the
singular extremal through $z$ and $\gamma$ the correspondent singular
curve. We set:
$$
J_\lambda(t)=D_\lambda F\left(\pi^{-1}\Delta_{\gamma(t)}\right)
\subset T_zC^0_\Delta/T_z\lambda .
$$
Then $J_\lambda(t)\supset J^0_\lambda(t)$ and $J_\lambda(t)$ is a
Lagrangian subspace of the symplectic space
$T_zC^0_\Delta/T_z\lambda$. In other words,
$J_\lambda(t)^\angle=J_\lambda(t)$, where
$$
\mathcal S^\angle\stackrel{\mathrm{def}}{=}\{\zeta\in
T_zC^0_\Delta: \sigma(\zeta,\mathcal S)=0\},\ \mathcal S\subset
T_z.
$$

Given $s\in\mathbb R\setminus\{0\}$, $s\lambda$ is the singular
extremal through $sz\in C_\Delta^0$. Hence $T_z(\mathbb Rz)\subset
J_\lambda(t),\ \forall t$ and $J_\lambda(t)\subset T_z(\mathbb
Rz)^\angle$. This inclusion allows to make one more useful
reduction. Namely, we set $\Sigma_z=T_z(\mathbb Rz)^\angle/ T_z\mathbb Rz$,
a symplectic space, $\dim\Sigma_z=2(n-3)$ and $J_\lambda(t)$
is a Lagrangian subspace of $\Sigma_z$.

Let $L(\Sigma_z)$ be the Lagrange Grassmannian: the manifold formed by
all Lagrangian subspaces of $\Sigma_z$. The curve
$t\mapsto J_\lambda(t)$ considered as a curve in $L(\Sigma_z)$
is called the {\it Jacobi curve} associated to the extremal $\lambda$.

The dimension of a Lagrangian subspace is one half of the dimension
of the ambient symplectic space. In particular, generic pair of
Lagrangian subspaces has zero intersection. Jacobi curves
are not at all generic, nevertheless, under very mild regularity
assumption on the distribution (see \cite{ze})  they satisfy the following
important property: $J_\lambda(t)\cap J_\lambda(\tau)=0$ for
sufficiently small $|t-\tau|\ne 0$.

Let $\pi_{t\tau}$ be the linear projector of $\Sigma_z$ on
$J_\lambda(\tau)$ along $J_\lambda(t)$. In other words,
$\pi_{t\tau}:\Sigma_z\to\Sigma_z$,
$$
\pi_{t\tau}\bigr|_{J_\lambda(t)}=0,\quad
\pi_{t\tau}\bigr|_{J_\lambda(\tau)}=\mathbf 1.
$$
\begin{lemma}[see\cite{az}] We have:
$$
\mathrm{tr}\left(\frac{\partial^2\pi_{t\tau}}{\partial
t\partial\tau}\right)=
\frac{(n-3)^2}{(t-\tau)^2}+g_\lambda(t,\tau),
$$
where $g_\lambda(t,\tau)$ is a symmetric function of $(t,\tau)$
that is smooth in a neighborhood of $(t,t)$, for all $t$
out of a discrete subset of the domain of $J_\lambda(\cdot)$.
\end{lemma}

In what follows, we tacitly assume that the value of $t$ is taken
out of the discrete subset provided by the lemma. A basic
invariant of the parameterized singular extremal $\lambda$ is the
{\it generalized Ricci curvature} $$\mathfrak
r_\lambda(\lambda(t)) \stackrel{\mathrm{def}}{=}g_\lambda(t,t).$$

Generalized Ricci curvature depends on the parameterization of the
extremal; this dependence is controlled by the following chain
rule. Let \linebreak $\varphi:\mathbb R\to\mathbb R$ be a change
of the parameter, then:
$$\mathfrak r_{\lambda\circ\varphi}\left(\lambda(\varphi(t))\right)=
\mathfrak
r_\lambda\left(\lambda(\varphi(t))\right)\dot\varphi^2(t)+(n-3)^2
{\mathbb S}(\varphi),
$$
where ${\mathbb S}(\varphi)=
\frac{\stackrel{\ldots}{\varphi}(t)}{2\dot\varphi(t)}-\frac 34
\left(\frac{\ddot\varphi(t)}{\dot\varphi(t)}\right)^2,$ the Schwartzian
derivative. The chain rule implies that the generalized
Ricci curvature $\mathfrak r_\lambda$ can be always made zero by a local
reparameterization of the extremal $\lambda$.
We say that a local parameter $t$ is a {\it projective parameter} if
$\mathfrak r_\lambda(t)\equiv 0$; such a parameter
is defined up to a M\"obius transformation.

Let $t$ be a projective parameter, then the quantity:
$$
A(\lambda(t))=\frac{\partial^2g}{\partial\tau^2}(t,\tau)\Bigr|_{\tau=t}
(dt)^4
$$
is a well defined degree 4 differential on $\lambda$; we
called it the {\it fundamental form} on $\lambda$.

In arbitrary, not necessary projective parameter, the fundamental
form has the following expression:
$$
A(\lambda(t))=
\left(\frac{\partial^2g}{\partial\tau^2}\Bigr|_{\tau=t}- \frac
3{5(n-3)^2}\mathfrak r_\lambda(t)^2-\frac 32\ddot\mathfrak r_\lambda(t)\right)
(dt)^4.
$$
Assume that $A(\lambda(t))\ne 0$, then the identity
$|A(\lambda(s))\left(\frac d{ds}\right)|=1$ defines a unique (up
to a translation) {\it normal parameter} $s$.

Let $z\in C^0_\Delta$ and $\lambda_s$ be the normally
parameterized singular extremal through $z$. We set
$$\bar\mathfrak r(z)=\mathfrak r_{\lambda_s}(z),$$ the {\it
projective  generalized Ricci curvature}.
Then $z\mapsto\bar\mathfrak r(z)$ is a
function on $C^0_\Delta$ which depends only on $\Delta$.

\medskip
Now come back to the case $k=2,\ n=5$. In this case, the
fundamental form $A$ is reduced to the famous Cartan's degree four
form on the distribution constructed in \cite{ca} by the method of
equivalence. The distribution is flat if and only if $A\equiv 0$.

Zelenko \cite{ze} performed detailed calculations for the rolling
balls model. As before, let $\rho$ be the ratio of the radii of
the balls. We assume that $1<\rho\le+\infty$. It appears that
$$\mathrm{sign}(A)=\mathrm{sign}(\rho-3).$$ Singular curves are just rolling
motions along geodesics (i.e. along big circles). The symmetry group acts
transitively on the space of geodesics, hence function $\mathfrak r$ must be
constant in this case. We have:
$$
\bar\mathfrak r=\frac{4\sqrt{35}(\rho^2+1)}{3\sqrt{(\rho^2-9)(9\rho^2-1)}}.
$$
In particular, the distributions corresponding to different $\rho$
are mutually non equivalent and only the distribution
corresponding to $\rho=3$ is flat.


\begin{thebibliography}{99}

\bibitem{aII} A.~Agrachev, {\it Feedback--invariant optimal control
theory and differential geometry, II. Jacobi curves for singular extremals.}
J. Dynamical and Control Systems, 1998, v.4, 583--604

\bibitem{as} A. Agrachev, Yu. Sachkov, {\it Control theory from the geometric
viewpoint}. Springer Verlag, 2004, xiv+412 pp.

\bibitem{az} A.~Agrachev, I.~Zelenko, {\it Geometry of Jacobi
curves, I, II}. J. Dynamical and Control Systems, 2002, v.8,
93--140; 167--215

\bibitem{aznu} A. Agrachev, I. Zelenko, {\it Nurovski's conformal structure
for (2,5)-distributions via dynamics of abnormal extremals}. Proceed.
RIMS workshop"Developments of Cartan geometry and related mathematical problems";
Kyoto, 2006, 204--218

\bibitem{br} R. Bryant, L. Hsu, {\it Rigidity of integral curves of rank 2
distributions}. Invent. Math., 1993, v.114, 435--461

\bibitem{ca} E. Cartan, {\it Les syst\`emes de Pfaff \`a cinq variables et les
\'equations aux d\'eriv\'ees partielles du second ordre}. Ann Sci. Ecole
Normale, 1910, v.27(3), 109--192

\bibitem{ga} R. Gardner, {\it The method of equivalence and its applications}.
SIAM, 1989 viii+127 pp.

\bibitem{ju} V. Jurdjevic, {\it Geometric control theory}. Cambridge Univ.
Press, 1997, xviii+492 pp.

\bibitem{ka} A. Kaplan, F. Levstein, {\it A split Fano plane}. In
preparation

\bibitem{mo} R. Mongomery, {\it A tour of subRiemannian geometries, their
geodesics and applications}. Amer. Math Soc., 2002, xx+259 pp.

\bibitem{nu} P. Nurowski, {\it Differential equations and conformal structures}.
J. Geometry and Physics, 2005, v.55, 19--49

\bibitem{ks} K. Sagerschnig, {\it Split octonions and generic rank 2 distributions in dimension 5}.
Proceed. 26th Winter School on Geometry and
Physics, Srni 2006" Rend. Crc. Mat. Palermo Suppl. II, to appear

\bibitem{sp} T. Springer, F. Veldkamp, {\it Octonians, Jordan algebras and
exceptional groups}. Springer Verlag, 2000, viii+208 pp.

\bibitem{ze} I.~Zelenko, {\it Variational approach to differential invariants
of rank 2 vector distributions}. J. Differential Geometry and Appl.,
2006, v.24, 235--259

\end{thebibliography}
\end{document}